\documentclass{amsart}
\usepackage{amsfonts}
\usepackage{amsmath}
\usepackage{amssymb}

\theoremstyle{plain}

\newtheorem{remark}{Remark}

\numberwithin{equation}{section}

\input{tcilatex}

\begin{document}
\title[Riemann zeta function and Wirtinger-type Inequalities]{\textbf{Large
Spaces Between the Zeros of the Riemann Zeta-Function }}
\author{\textbf{S. H. Saker}}
\address{{\small Department of Mathematics Skills, PYD, King Saud
University, Riyadh 11451, } {\small Saudi Arabia, Department of Math.,
Faculty of Science, Mansoura University, Mansoura 35516, Egypt.}}
\email{{\small shsaker@mans.edu.eg, mathcoo@py.ksu.edu.sa}}
\subjclass{11M06, 11M26.}
\keywords{Riemann zeta function, zeros the Riemann zeta function }
\maketitle

\begin{abstract}
In this paper, we will employ the Opial and Wirtinger type inequalities to
derive some conditional and unconditional lower bounds for the gaps between
the zeros of the Riemann zeta-function. First, we prove\ (unconditionally)
that the consecutive nontrivial zeros often differ by at least $1.9902$
times the average spacing. This value improves the value $1.9$ due to
Mueller and the value $1.9799$ due to Montogomery and Odlyzko. Second, on
the hypothesis that the $2k-$th mixed moments of the Hardy $Z-$function and
its derivative are correctly predicted by random matrix theory, we derive
some explicit formulae for the gaps and use them to establish new
(conditional) large gaps.
\end{abstract}

\section{Introduction}

The Riemann zeta function $\zeta (s)$ is defined on $\{s\in \mathbb{C}:\func{%
Re}(s)>1\}$ by the series\ 
\begin{equation*}
\zeta (s):=1+\frac{1}{2^{s}}+\frac{1}{3^{s}}+\frac{1}{4^{s}}+...,\text{ \ \
\ \ \ for \ }\func{Re}s>1,
\end{equation*}%
which converges in the region described by the Cauchy integral test. There
is another representation of $\zeta $ due to Euler in 1749 which is perhaps
more fundamental and which is the reason for the significance of the zeta
function and gives analytic expression to the fundamental theorem of
arithmetic. This formula is given by 
\begin{equation*}
\zeta (s):=\dprod_{p}\left( 1-\frac{1}{p^{s}}\right) ^{-1},\ \ \ \text{for\ }%
\ \ \ \func{Re}s>1,
\end{equation*}%
where the product is taken over all prime numbers. The zeta-function is one
of the most studied transcendental functions, having in view its many
applications in number theory, algebra, complex analysis, statistics as well
as in physics. Another reason why this function has drawn so much attention
is the celebrated Riemann conjecture regarding nontrivial zeros which states
that all nontrivial zeros of the Riemann zeta function $\zeta (s)$ lie on
the critical line $\func{Re}(s)=1/2.$ Riemann showed that the zeta-function
satisfies a functional equation of the form%
\begin{equation}
\pi ^{-s/2}\Gamma (\frac{s}{2})\zeta (s)=\pi ^{-(1-s)/2}\Gamma (\frac{1-s}{2}%
)\zeta (1-s),  \label{F}
\end{equation}%
where $\Gamma $ is the Euler gamma function. By this equation there exist-so
called trivial (real) zeros at $s=-2n$ for any positive integer $n$
(corresponding to the poles of the appearing Gamma-factors), and all
nontrivial (non-real) zeros are distributed symmetrically with respect to
the critical line $\func{Re}s=1/2$ and the real axis. The significance
contribution of the formula (\ref{F}) is the consideration of the
zeta-function as an analytic function. We note from the functional equation
that if $\rho \in \mathbf{Z}$ is a zero of $\zeta (s)$, so is $1-\rho $, $%
\overline{\rho }$, $1-\overline{\rho }$ and according to the Riemann
hypothesis $\func{Re}\rho =1/2$ and under this hypothesis $1-\rho =\overline{%
\rho }$ and $\zeta (\overline{\rho })=\overline{\zeta (\rho )}$. Clearly,
there are no zeros in the half-plane of convergence $\func{Re}(s)>1,$ and it
is also known that $\zeta (s)$ does not vanish on the line $\func{Re}(s)=1.$

The number $N(t)$ of the non-trivial zeros of $\zeta (s)$ with ordinate in
the interval $[0,$ $T]$ is asymptotically given by the Riemann-von Mangoldt
formula (see \cite{RS}) 
\begin{equation*}
N(T)=\frac{T}{2\pi }\log (\frac{T}{2\pi e})+O(\log T).
\end{equation*}%
Consequently there are infinitely many nontrivial zeros, all of them lying
in the critical strip $0<\func{Re}s<1,$ and the frequency of their
appearance is increasing as $T\rightarrow \infty .$

There are three directions regarding the studies of the zeros of the Riemann
zeta function. The first direction is concerning with the existence of
simple zeros. It is conjectured that all or at least almost all zeros of the
zeta-function are simple. For this direction Conrey \cite{C2} proved that
more than two-fifths of the zeros are simple and on the critical line. This
value has been improved by Cheer and Goldston \cite{CG} and proved that at
least $0.662753$ of the zeros are simple assuming the Riemann hypothesis.
The second direction is the most important goal, is the determination of the
moments of the Riemann zeta function on the critical line and the evaluation
of the Riemann zeta function at integers which gives an integeral
representaion of this function. It is important because it can be used to
estimate the maximal order of the zeta-function on the critical line, and
because of its applicability in studying the distribution of prime numbers
and to divisor problems. For more details of this direction, we refer the
reader to \cite{JCAM}, \cite{Hughes} and \cite{S} and the references cited
therein. The third direction is the distribution of the zeros when the
Riemann hypothesis is satisfied which is of our interest in this paper. In
fact the distribution of zeros of the Riemann zeta-function is of
fundamental importance in number theory as well as in physics. In the
following, we briefly present some results related to this direction. Assume
that $(\beta _{n}+i\gamma _{n})$ are the zeros of $\zeta (s)$ in the upper
half-plane (arranged in nondecreasing order and counted according to
multiplicity) and $\gamma _{n}\leq \gamma _{n+1}$ are consecutive ordinates
of all zeros. We put 
\begin{equation}
r_{n}:=\frac{(\gamma _{n+1}-\gamma _{n})}{(2\pi /\log \gamma _{n})},
\label{Q}
\end{equation}%
and define $\lambda :=\lim \sup_{n\rightarrow \infty }r_{n},\ $and\ $\mu
:=\lim \inf_{n\rightarrow \infty }r_{n}.$ The numbers $\mu $ and $\lambda $
have received a great deal of attention. As mentioned by Montogomery \cite%
{M1} it would be interesting to see how numerical evidence compare with the
above conjectures. It generally conjectured that 
\begin{equation}
\mu =0,\text{ \ \ and \ \ \ }\lambda =\infty .  \label{b}
\end{equation}%
Now, several results has been obtained, however the failure of Gram's low
(see \cite{H}) indicate that the asymptotic behavior is approached very
slowly. Thus the numerical evidence may not be particularly illuminating. In
fact, important results concerning the values of $\lambda $ and $\gamma $
have been obtained by some authors. Selberg \cite{Selberg} proved that $%
0<\mu <1<\lambda ,$ and the average of $r_{n}$ is $1.$ Note that $2\pi /\log
\gamma _{n}$ is the average spacing between zeros. Fujii \cite{Fujii} also
showed that there exist constants $\lambda >1$ and $\mu <1$ such that 
\begin{equation*}
\frac{(\gamma _{n+1}-\gamma _{n})}{(2\pi /\log \gamma _{n})}\geq \lambda 
\text{, and \ }\frac{(\gamma _{n+1}-\gamma _{n})}{(2\pi /\log \gamma _{n})}%
\leq \mu ,
\end{equation*}%
each holds for a positive proportion of $n.$ Mueller \cite{33} obtained $%
\lambda >1.9,$ assuming the Riemann hypothesis. Montogomery and Odlyzko \cite%
{Montgomery} showed, assuming the Riemann hypothesis, that $\lambda >1.9799,$
and $\mu <0.5179.$ Conrey, Ghosh and Gonek \cite{Conrey1} proved that if the
Riemann hypothesis is true, then $\lambda >2.337,$ and\ $\mu <0.5172.$
Conrey, Ghosh and Gonek \cite{Conrey2} obtained a new lower bound and proved
that $\lambda >2.68,$ assuming the generalized Riemann hypothesis for the
zeros of the Dirichlet $L-$ functions. Bui, Milinovich and Ng \cite{BMN}
obtained $\lambda >2.69,$ and $\mu <0.5155,$ assuming the Riemann
hypothesis. Ng in \cite{Ng} proved that $\lambda >3,$ assuming the
generalized Riemann hypothesis for the zeros of the Dirichlet $L-$functions.
Note that any other small values of $\mu $ and large values of $\lambda $
will help in proving the conjecture (\ref{b}).

Let $\Lambda $ denote the quantity in (\ref{Q}) where only zeros $\frac{1}{2}%
+it_{n}$ on the critical line, i.e., we define 
\begin{equation}
\Lambda :=\lim \sup \frac{t_{n+1}-t_{n}}{(2\pi /\log t_{n})}.  \label{R2}
\end{equation}%
Note that the Riemann hypothesis implies that the $t_{n}$ corresponded to
the positive ordinates of non-trivial zeros of the zeta function, i.e., $%
N(T)\sim \frac{T}{2\pi }\log T.$ The average spacing between consecutive
zeros with ordinates of order $T$ is $2\pi /\log (T)$ which tends to zero as 
$T\rightarrow \infty .$ Hall \cite{Hall3} showed that $\Lambda \geq \lambda $%
, and the lower bound for $\Lambda $ bear direct comparison with such bounds
for $\lambda $ dependent on the Riemann hypothesis, since if this were true
the distinction between $\Lambda $ and $\lambda $ would be nugatory. Of
course $\Lambda \geq \lambda $ and the equality holds if the Riemann
hypothesis is true. The behavior of $\zeta (s)$ on the critical line is
reflected by the Hardy $Z-$function $Z(t)$ as a function of a real variable,
defined by 
\begin{equation}
Z(t)=e^{i\theta (t)}\zeta (\frac{1}{2}+it),\text{ where }\theta (t):=\pi
^{-it/2}\frac{\Gamma (\frac{1}{4}+\frac{1}{2}it)}{\left\vert \Gamma (\frac{1%
}{4}+\frac{1}{2}it)\right\vert }.  \label{Z}
\end{equation}%
It follows from the functional equation (\ref{F}) for $\zeta (s)$ that $Z(t)$
is an infinitely often differentiable and real function for real $t.$
Moreover $\left\vert Z(t)\right\vert =\left\vert \zeta (1/2+it)\right\vert $%
. Consequently, the zeros of $Z(t)$ correspond to the zeros of the Riemann
zeta-function on the critical line. The moments $I_{k}(T)$ of the Hardy $Z-$%
function $Z(t)$ function and the moments $M_{k}(T)$ of its derivative are
defined by 
\begin{equation*}
I_{k}(T):=\int_{0}^{T}\left\vert Z(t)\right\vert ^{2k}dt,\text{ and }%
M_{k}(T):=\int_{0}^{T}\left\vert Z^{^{\prime }}(t)\right\vert ^{2k}dt.
\end{equation*}%
For positive real numbers $k$, it is believed that $I_{k}(T)\sim C(k)$ $%
T\left( \log T\right) ^{k^{2}}$and $M_{k}(T)\sim L(k)T\left( \log T\right)
^{k^{2}+2k}$ for positive constants $C_{k}$ and $L_{k}$ will be defined
later. Keating and Snaith \cite{KS} based on considerations from random
matrix theory conjectured that%
\begin{equation}
I_{k}(T)\sim a(k)b(k)T\left( \log T\right) ^{k^{2}},  \label{A1}
\end{equation}%
where $a(k)$ and $b(k)$ are defined by%
\begin{equation}
a(k):=\tprod_{p}(\left( 1-\frac{1}{p^{2}}\right) \dsum_{m=0}^{\infty }\left( 
\frac{\Gamma (m+k)}{m!\Gamma (k)}\right) ^{2}p^{-m},\text{ and }%
b(k):=\dprod_{j=0}^{k-1}\frac{j!}{(j+k)!}.  \label{B1}
\end{equation}%
To find the lower bound of $\Lambda $ Hall \cite{Hall1} used a
Wirtinger-type inequality of Beesack \cite{B} and the moment 
\begin{equation}
\int_{0}^{T}Z^{4}(t)dt=\frac{1}{2\pi ^{2}}T\log ^{4}(t)+O(T\log ^{3})\text{, 
}  \label{Hl1}
\end{equation}%
due to Ingham (\cite{Ingham}) and the moment%
\begin{equation}
\int_{0}^{T}(Z^{^{\prime }}(t))^{4}dt=\frac{1}{1120\pi ^{2}}T\log
^{8}(t)+O(T\log ^{7}),  \label{Hl2}
\end{equation}%
due to Conrey \cite{Conrey3} and proved unconditionally that 
\begin{equation}
\Lambda \geq (\frac{105}{4})^{\frac{1}{4}}=2.2635.  \label{D0}
\end{equation}%
In \cite{Hall2} Hall remarked that Beesack inequality is sharp but it is not
optimal for application and proved a new Wirtinger-type inequality and used
the moments (\ref{Hl1})-(\ref{Hl2}) and the moment 
\begin{equation}
\int_{0}^{T}Z^{2}(t)(Z^{^{\prime }}(t))^{2}dt=\frac{1}{120\pi ^{2}}T\log
^{6}(t)+O(T\log ^{5}),  \label{Hl3}
\end{equation}%
due to Conrey \cite{Conrey3}, and proved (unconditionally) that $\Lambda
\geq \sqrt{11/2}=2.345\,2.$ The moments of $Z(t)$ and its derivative (of
mixed powers)%
\begin{equation}
\int_{0}^{T}Z^{2k-2h}(t)(Z^{^{\prime }}(t))^{2h}dt\sim C(h,k)T\left( \log
T\right) ^{k^{2}+2h},  \label{M1}
\end{equation}%
has been predicted by Random Matrix Theory (RMT) by Hughes \cite{Hughes} who
stated an interesting conjecture on the moments subject to the truth of
Riemann's hypothesis when the zeros are simple. This conjecture includes for
fixed $k>-3/2$ the asymptotes formula of the moments of the higher order of
the Riemann zeta function and its derivative. We suppose further that if $k$
is a fixed positive integer and $h\in \lbrack 0,$ $k]$ is an integer then
the formula 
\begin{equation}
\int_{0}^{T}Z^{2k-2h}(t)(Z^{^{\prime }}(t))^{2h}dt\sim a(k)b(h,k)T\left(
\log T\right) ^{k^{2}+2h},  \label{Hu}
\end{equation}%
holds. Note that (\ref{Hu}) has benn predicted by Keating and Snaith \cite%
{KS} in the case when $h=0$, with wider range $\func{Re}(k)>-1/2$ and by
Hughes \cite{Hughes} in the range $\min (h,k-h)>-1/2,$ $a(k)$ is a product
over the primes and $b(h,k)$ is rational: indeed for integral $h$, it is
obtained that 
\begin{equation}
b(h,k):=b(0,k)\left( \frac{\left( 2h\right) !}{8^{h}h!}\right) H(h,k),
\label{BG}
\end{equation}%
where $H(h,k)$ is an explicit rational function of $k$ for each fixed $h$
and $b(0,k)=b(k)$ which is defined as in (\ref{B1}). The functions $H(h,k)$
as introduced by Hughes \cite{Hughes} are given in the following table where 
$K=2k$:

\begin{center}
\begin{equation*}
\begin{tabular}[t]{|l|}
\hline
$H(0,k)=1,$ $H(1,k)=\frac{1}{K^{2}-1},$ $H(2,k)=\frac{1}{(K^{2}-1)(K^{2}-9)}%
, $ \\ \hline
$H(3,k)=\frac{1}{(K^{2}-1)^{2}(K^{2}-25)},$ $H(4,k)=\frac{K^{2}-33}{%
(K^{2}-1)^{2}(K^{2}-9)(K^{2}-25)(K^{2}-49)},$ \\ \hline
$H(5,k)=\frac{K^{4}-90K^{2}+1497}{%
(K^{2}-1)^{2}(K^{2}-9)^{2}(K^{2}-25)(K^{2}-49)(K^{2}-81)},$ \\ \hline
$H(6,k)=\frac{K^{6}-171K^{4}+6867K^{2}-27177}{%
(K^{2}-1)^{3}(K^{2}-9)^{2}(K^{2}-25)(K^{2}-49)(K^{2}-81)(K^{2}-121)},$ \\ 
\hline
$H(7,k)=\frac{K^{8}-316K^{6}+30702K^{4}-982572K^{2}+6973305}{%
(K^{2}-1)^{3}(K^{2}-9)^{2}(K^{2}-25)^{2}(K^{2}-49)(K^{2}-81)(K^{2}-121)(K^{2}-169)%
},$ \\ \hline
\end{tabular}%
\end{equation*}%
\textit{Table 1. The values of }$H(h,k)$\textit{, for }$h=0,1,2,...7$ 
\textit{where }$K=2k.$
\end{center}

This sequence continuous, and it is believed that both the nominator and
denominator are polynomials in $k^{2}$, moreover that the denominator is
actually (see \cite{D}) 
\begin{equation}
\dprod_{a\text{ }odd>0}\left\{ (K^{2}-a^{2})^{\alpha (a,h)}:\alpha (a,h)=%
\frac{4h}{a+\sqrt{a^{2}+8h}}\right\} .  \label{Monic}
\end{equation}%
Using the equation (\ref{BG}) and the definitions of the functions $H(h,k)$,
we can obtain the values of $b(0,k)/b(k,k)$ for $k=1,2,...,7.$ Hall \cite%
{Hall8} shown that in the case when $h=3$, ($H(3,k))$ requires adjustment to
fit with (\ref{Monic}) in that extra factor $K^{2}-9$ should be introduced
in both the nominator and denominator. I hope also to get the other values
of $H(h,k)$ for $k\geq 8$ which will help in deriving new values of $\Lambda 
$, since our calculation (at this moment) will stop at $k=7.$ The values of $%
b(0,k)/b(k,k)$ for $\ k=1,2,...,7,$ that we will need in this paper are
determined from (\ref{BG}) and presented in the following table:

\begin{center}
\begin{equation*}
\begin{tabular}[t]{|l|l|l|l|}
\hline
$\frac{b(0,1)}{b(1,1)}$ & $\frac{b(0,2)}{b(2,2)}$ & $\frac{b(0,3)}{b(3,3)}$
& $\frac{b(0,4)}{b(4,4)}$ \\ \hline
$12$ & $\frac{6720}{12}$ & $\frac{49674240}{864}$ & $\frac{271159356948480}{%
31(870912)}$ \\ \hline
$\frac{b(0,5)}{b(5,5)}$ & $\frac{b(0,6)}{b(6,6)}$ & $\frac{b(0,7)}{b(7,7)}$
&  \\ \hline
$\frac{581050229760}{227}$ & $\frac{114664452340838400}{133933}$ & $\frac{%
1769682901766011323008}{5078125}$ &  \\ \hline
\end{tabular}%
\end{equation*}

\textit{Table 2. The Values of the }$b(0,k)/b(k,k)$\textit{\ for }$\
k=1,2,...,7.$
\end{center}

Hall in \cite{Hall3, Hall8} used the moments of mixed powers (\ref{Hu}) and
a new Wirtinger-type inequality designed exclusively for this problem to
improve the lower values of $\Lambda $. In particular Hall \cite{Hall3}
proved a the Wirtinger-type inequality 
\begin{equation}
\int_{0}^{\pi }H\left( y^{^{\prime }}(t)/y(t)\right) y^{2k}(t)dt\geq
(2k-1)L\int_{0}^{\pi }y^{2k}(t)dt,  \label{W2}
\end{equation}%
where $L=L(k,H)$ is determined from the solution of the equation 
\begin{equation*}
\int_{0}^{\infty }\frac{G^{^{\prime }}(u)}{G(u)+(2k-1)L}\frac{du}{u}=k\pi 
\text{, for }k\in \mathbb{N},
\end{equation*}%
where $G(u):=uH^{^{\prime }}(u)-H(u)$, $y=y(t)\in C^{2}[0,$ $\pi ]$ and $%
y(0)=y(\pi )=0$, $H(u)$ be an even function, increasing, strictly convex on $%
\mathbb{R}^{+}$ and satisfies $H(0)=H^{^{\prime }}(0)=0\ $and $uH^{^{\prime
\prime }}(u)\rightarrow 0$ as $u\rightarrow 0.$ The inequality (\ref{W2}) is
proved by using the calculus of variation which depends on the minimization
of the integral on the left hand side subject to the constrains $y(0)=0$ and 
$\int_{0}^{\pi }y^{2k}(t)dt=1.$ Assuming that\ (\ref{Hu}) is correctly
predicted, Hall employed the inequality (\ref{W2}) when 
\begin{equation*}
H(u):=\dsum_{h=1}^{k}\frac{2k-1}{2h-1}\left( 
\begin{array}{c}
h \\ 
k%
\end{array}%
\right) \upsilon _{h}u^{2h}\text{, \ \ }\upsilon _{h}\geq 0\text{, \ \ }%
\upsilon _{k}=1,
\end{equation*}%
and obtained 
\begin{equation}
\Lambda \geq \sqrt{7533/901}=2.8915.  \label{D2}
\end{equation}%
The main challenge in \cite{Hall3} was to maximize $X=\kappa ^{2}$ (which is
not an easy task) where $X$ satisfies the equation $27X^{3}+385\mu
X^{2}+10395\vartheta X-121275L=0,$ and $L$ obtained form the equation 
\begin{equation*}
\int_{-\infty }^{\infty }\frac{x^{4}+2\mu x^{2}+\upsilon }{x^{6}+3\mu
x^{4}+3\upsilon x^{2}+L}dx=\pi .
\end{equation*}%
In \cite{Hall4} Hall employed the generalized Wirtinger inequality (\ref{W2}%
) and simplified the calculation in \cite{Hall3} and converted the problem
into one in the classical theory of equations involving Jacobi-Schur
functions to maximize $X$. Assuming that (\ref{Hu}) is correctly predicted
by RMT, Hall obtained the new values of $\Lambda $ which is listed in the
following table:%
\begin{equation}
\begin{tabular}[t]{|l|l|l|l|}
\hline
$\Lambda (3)$ & $\Lambda (4)$ & $\Lambda (5)$ & $\Lambda (6)$ \\ \hline
$\sqrt{7533/901}$ & $3.392272$ & $3.858851$ & $4.2981467.$ \\ \hline
\end{tabular}
\label{Hall}
\end{equation}%
The methods that have been used by Hall to establish the lower bounds of $%
\Lambda $ are quite complicated and need a lot of calculations as well as
the reader should be familiar with calculus of variations and optimization
theory. In \cite{R} the authors applied a technique involving the comparison
of the continuous global average with local average obtained from the
discrete average to a problem of gaps between the zeros of zeta function
assuming the Riemann hypothesis. Using this approach, which takes only zeros
on the critical line into account, the authors computed similar bounds under
assumption of the Riemann hypothesis when (\ref{Hu}) holds. In particular
they showed that for fixed positive integer $r$%
\begin{equation}
\frac{(\gamma _{n+r}-\gamma _{n})}{(2\pi r/\log \gamma _{n})}\geq \theta ,
\label{Steuding}
\end{equation}%
holds for any $\theta \leq 4k/\pi re$ for more than $c(\log T)^{-4k^{2}}$
proportion of the zeros $\gamma _{n}\in \lbrack 0,T]$ with a computable
constant $c=c(k,\theta ,r).$

\bigskip

In this paper, we\ will employ some well-known Opial and Wirtinger type
inequalities to derive new unconditional lower bound for $\Lambda $ and also
establish some explicit formulae for the gaps between the zeros. First, we
apply the Wirtinger type inequality due to Brneti\'{c} and Pe\v{c}ari\'{c} 
\cite{BP} and prove that $\Lambda \geq 1.9902$\ (unconditionally) which
improves the value $1.9$ of Mueller and the value $1.9799$ of Montogomery
and Odlyzko. Second, assuming that the moments of $Z(t)$ and its derivative
are correctly predicted by RMT, we established some new explicit formulae
for $\Lambda (k)$ by employing an Opial inequality due to Yang \cite{Yang}
and Wirtinger type inequality due to Agarwal and Pang \cite{AP}. As an
application, we derived some new conditional series of the lower bounds for $%
\Lambda (k)$. Our results do not require any additional information from the
calculus of variation and optimization theory.

\section{Main Results}

In this section, we\ employ some well-known Opial and Wirtinger type
inequalities to prove the main results. First, we employ the Wirtinger type
inequality due to Brneti\'{c} and Pe\v{c}ari\'{c} \cite{BP} to find a new
unconditional lower bound for $\Lambda .$ The Wirtinger type inequality due
to Brneti\'{c} and Pe\v{c}ari\'{c} \cite{BP} is presented in the following
theorem.

\bigskip

\textbf{Theorem A}. \textit{Assume that }$x(t)\in C^{1}[0,$\textit{\ }$\pi ]$%
\textit{\ and }$x(0)=x(\pi )=0$\textit{, then}%
\begin{equation}
\int_{0}^{\pi }(x^{^{\prime }}(t))^{2k}dt\geq \frac{1}{\pi ^{2k}I(k)}%
\int_{0}^{\pi }x^{2k}(t)dt,\text{ \ for \ }k\geq 1,  \label{AP2}
\end{equation}%
\textit{where }%
\begin{equation*}
I(k)=\int_{0}^{1}\frac{1}{(t^{1-2k}+(1-t)^{1-2k})}dt.
\end{equation*}

In the following, we will apply the inequality (\ref{AP2}) and the moment (%
\ref{Hl1}) due to Ingham \cite{Ingham}) and the moment (\ref{Hl2}) due to
Conrey \cite{Conrey3} to find the new unconditional value of $\Lambda .$
From (\ref{AP2}), when $k=2,$ we have 
\begin{equation}
\int_{0}^{\pi }(x^{^{\prime }}(t))^{4}dt\geq \frac{1}{\pi ^{4}I(2)}%
\int_{0}^{\pi }x^{4}(t)dt,\   \label{AP3}
\end{equation}%
where 
\begin{equation*}
I(2):=\int_{0}^{1}\frac{1}{(t^{1-4}+(1-t)^{1-4})}dt=\frac{2863}{125000}.
\end{equation*}%
By a suitable linear transformation, we can deduce from (\ref{AP3}) that if $%
x(t)\in C^{1}[a,b]$ and $x(a)=x(b)=0,$ then 
\begin{equation}
\int_{a}^{b}(\frac{b-a}{\pi })^{4}(x^{^{\prime }}(t))^{4}dt\geq \frac{125000%
}{2863\pi ^{4}}\int_{a}^{b}x^{2k}(t)dt  \label{2.3}
\end{equation}%
where $x(a)=x(b)=0.$

\bigskip

The following theorem gives the new unconditional value of $\Lambda .$

\bigskip

\textbf{Theorem 2.1.} \textit{Let }$\varepsilon (T)\rightarrow 0$\textit{\
in such a way that }$\varepsilon (T)\log T\rightarrow \infty .$\textit{\
Then for sufficiently large }$T$\textit{, there exists an interval contained
in }$[T,(1+\varepsilon (T))T]$\textit{\ which is free of zeros of }$Z(t)$%
\textit{\ and having length at least }%
\begin{equation*}
\frac{1}{2\pi }\sqrt[4]{\frac{10000000}{409}}\left\{ 1+O\left( \frac{1}{%
\varepsilon (T)\log T}\right) \right\} \frac{2\pi }{\log T}.
\end{equation*}%
\textit{Thus }%
\begin{equation}
\Lambda \geq 1.9902.  \label{new}
\end{equation}%
\textbf{Proof.} We follow the arguments in \cite{Hall2} to prove our
theorem. Suppose that $t_{l}$ is the first zero of $Z(t)$ not less than $T$
and $t_{m}$ the last zero not greater than $(1+\varepsilon )T$ where $%
\varepsilon (T)\rightarrow 0$ in such a way that $\varepsilon (T)\log
T\rightarrow \infty .$ Suppose further that for $l\leq n<m$, we have 
\begin{equation}
L_{n}=t_{n+1}-t_{n}\leq \frac{2\pi \kappa }{\log T}.  \label{kv3}
\end{equation}%
Applying the inequality (\ref{2.3}) with $a=t_{n}$, $b=t_{n+1}$ and $%
y(t)=Z(t)$, we have 
\begin{equation*}
\int_{t_{n}}^{t_{n+1}}(\frac{L_{n}}{\pi })^{4}(Z^{^{\prime }}(t))^{4}dt-%
\frac{1}{\pi ^{4}}\frac{125\,000}{2863}\int_{a}^{b}Z^{4}(t)dt\geq 0.
\end{equation*}%
Since the inequality remains true if we replace $L_{n}/\pi $ by $2\kappa
/\log T$, we have 
\begin{equation}
\int_{t_{n}}^{t_{n+1}}\left[ \left( \frac{2\kappa }{\log T}\right)
^{4}\left( Z^{^{\prime }}(t)\right) ^{4}-\frac{1}{\pi ^{4}}\frac{125\,000}{%
2863}Z^{4}(t)\right] dt\geq 0.  \label{kv4}
\end{equation}%
Summing (\ref{kv4}) over $n,$ and using the moments (\ref{Hl1})-(\ref{Hl2}),
we obtain 
\begin{eqnarray*}
&&\frac{1}{1120\pi ^{2}}\left( \frac{2\kappa }{\log T}\right) ^{4}T\log
^{8}(T)+O(T\log ^{7}) \\
&&-\frac{1}{\pi ^{4}}\frac{125\,000}{2863}\frac{1}{2\pi ^{2}}T\log
^{4}(T)+O(T\log ^{3}T) \\
&=&\frac{\left( 2\kappa \right) ^{4}}{1120\pi ^{2}}T\log ^{4}(T)+O(T\log
^{3}T) \\
&&-\frac{1}{\pi ^{4}}\frac{125000}{2863}\frac{1}{2\pi ^{2}}(T\log
^{4}T)+O(T\log ^{3}T).
\end{eqnarray*}%
Follows the proof of Theorem 1 in \cite{Hall2}, we obtain 
\begin{equation*}
\kappa ^{4}\geq \frac{1120}{2^{5}\pi ^{4}I(2)}+O(1/\varepsilon (T)\log T).
\end{equation*}%
Then, we have (noting $(\varepsilon (T)\log T\rightarrow \infty $ as $%
T\rightarrow \infty )$ that 
\begin{equation*}
\Lambda \geq \frac{1}{2\pi }\sqrt[4]{\frac{125000}{2863}\frac{1120}{2}}%
=1.9902.
\end{equation*}%
The proof is complete.

\begin{remark}
One can easily see that the value $\Lambda \geq 1.9902$ improves the value $%
1.9$ of Mueller and the value $1.9799$ of Montogomery and Odlyzko.
\end{remark}

Next, in the following, we will apply the Opial inequality due to Yang \cite%
{Yang} to establish an explicit formula for the lower bounds of $\Lambda .$
The Yang inequality presented in the following theorem.

\bigskip

\textbf{Theorem B. }\textit{If }$x$\textit{\ is absolutely continuous on }$%
[a,$\textit{\ }$b]$\textit{\ with }$x(a)=0$\textit{\ (or }$x(b)=0)$\textit{,
then}%
\begin{equation}
\int_{a}^{b}\left\vert x(t)\right\vert ^{m}\left\vert x^{^{\prime
}}(t)\right\vert ^{n}dt\leq \frac{n}{m+n}(b-a)^{m}\int_{a}^{b}\left\vert
x^{^{\prime }}(t)\right\vert ^{m+n}dt,  \label{KV}
\end{equation}

\bigskip

The inequality (\ref{KV}) has immediate application to the case where $%
x(a)=x(b)=0$. Choose $c=(a+b)/2$ and apply (\ref{KV}) to $[a,c]$ and $[c,b]$
and then add to obtain%
\begin{eqnarray*}
&&\int_{a}^{b}\left\vert x(t)\right\vert ^{m}\left\vert x^{^{\prime
}}(t)\right\vert ^{n}dt \\
&\leq &\frac{n}{m+n}(\frac{b-a}{2})^{m}\left( \int_{a}^{c}\left\vert
x^{^{\prime }}(t)\right\vert ^{m+n}dt+\int_{c}^{b}\left\vert x^{^{\prime
}}(t)\right\vert ^{m+n}dt\right) \\
&\leq &\frac{n}{m+n}(\frac{b-a}{2})^{m}\left( \int_{a}^{b}\left\vert
x^{^{\prime }}(t)\right\vert ^{m+n}dt\right) .
\end{eqnarray*}%
So that if $x(0)=x(\pi )=0$, we have%
\begin{equation}
\int_{0}^{\pi }\left\vert x(t)\right\vert ^{m}\left\vert x^{^{\prime
}}(t)\right\vert ^{n}dt\leq \frac{n}{m+n}(\frac{\pi }{2})^{m}\int_{0}^{\pi
}\left\vert x^{^{\prime }}(t)\right\vert ^{m+n}dt.  \label{Kv3}
\end{equation}

\bigskip

\textbf{Theorem 2.2. }\textit{On the hypothesis that the Riemann hypothesis
is true and (\ref{Hu}) is correctly predicted, we have\ }%
\begin{equation}
\Lambda \geq \Lambda ^{\ast }(h,k)=\frac{1}{\pi }\left( \frac{k}{h}\frac{%
b(h,k)}{b(k,k)}\right) ^{\frac{1}{2k-2h}},\text{ for }h\neq k\neq 0.
\label{kli}
\end{equation}%
\textbf{Proof.} As in the proof of Theorem 2.2 by\textbf{\ }applying the
inequality (\ref{Kv3}) with $a=t_{n}$, $b=t_{n+1},$ $m=2k-2h$, $n=2h$, and $%
y=Z(t)$, we have 
\begin{eqnarray*}
&&\int_{t_{n}}^{t_{n+1}}\left( \frac{L_{n}}{\pi }\right) ^{2k}\left(
Z^{^{\prime }}(t)\right) ^{2k} \\
&\geq &\frac{k}{h}(\frac{2}{\pi })^{2k-2k}\int_{t_{n}}^{t_{n+1}}\left( \frac{%
L_{n}}{\pi }\right) ^{2h}\left\vert Z(t)\right\vert ^{2k-2h}\left\vert
Z^{^{\prime }}(t)\right\vert ^{2h}dt.
\end{eqnarray*}%
Since the inequality remains true if we replace $L_{n}/\pi $ by $2\kappa
/\log T$, we have%
\begin{eqnarray}
&&\int_{t_{n}}^{t_{n+1}}\left( \frac{2\kappa }{\log T}\right)
^{2k}\left\vert Z^{^{\prime }}(t)\right\vert ^{2k}  \notag \\
&\geq &\int_{t_{n}}^{t_{n+1}}\frac{k}{h}(\frac{2}{\pi })^{2k-2k}\left( \frac{%
2\kappa }{\log T}\right) ^{2h}\left\vert Z(t)\right\vert ^{2k-2h}\left\vert
Z^{^{\prime }}(t)\right\vert ^{2h}dt.  \label{BA3}
\end{eqnarray}%
Summing (\ref{BA3}) over $n,$ using (\ref{Hu}) we obtain%
\begin{eqnarray*}
&&\left( \frac{2\kappa }{\log T}\right) ^{2k}a(k)b(k,k)T\left( \log T\right)
^{k^{2}+2k} \\
&\geq &\frac{k}{h}(\frac{2}{\pi })^{2k-2k}\left( \frac{2\kappa }{\log T}%
\right) ^{2h}a(k)b(h,k)T\left( \log T\right) ^{k^{2}+2h}dt.
\end{eqnarray*}%
This implies that 
\begin{eqnarray*}
&&T\left( \log T\right) ^{k^{2}}\left\{ \left( 2\kappa \right)
^{2k}a(k)b(k,k)-\frac{k}{h}(\frac{2}{\pi })^{2k-2k}\left( 2\kappa \right)
^{2h}a(k)b(h,k)\right\} \\
&\geq &o(T\left( \log T\right) ^{k^{2}})
\end{eqnarray*}%
whence%
\begin{equation*}
\kappa ^{2k-2h}\geq \frac{1}{2^{2k-2h}}\frac{k}{h}(\frac{2}{\pi })^{2k-2k}%
\frac{b(h,k)}{b(k,k)}+o(1),\text{ (as }T\rightarrow \infty ).
\end{equation*}%
This implies that 
\begin{equation*}
\Lambda ^{2k-2h}(k)\geq \frac{1}{2^{2k-2h}}\frac{k}{h}(\frac{2}{\pi }%
)^{2k-2k}\frac{b(h,k)}{b(k,k)},\text{ }h\neq k\neq 0.
\end{equation*}%
which is the desired inequality and completes the proof.

\bigskip

To apply (\ref{kli}), we will need the following values of $b(1,k)$ and $%
b(k,k)$ that are determined from (\ref{BG}) where $H(h,k)$ are defined as in
Table 1.

\begin{eqnarray*}
b(1,2) &=&\tfrac{1}{720},\text{ \ }b(2,2)=\tfrac{1}{6720}\text{, }b(1,3)=%
\tfrac{1}{1209\,600},\text{ }b(3,3)=\tfrac{1}{496742400}, \\
b(1,4) &=&\tfrac{1}{219469\,824\,000}\text{, \ \ }b(4,4)=\tfrac{31}{%
271159356948480000}, \\
b(1,5) &=&\tfrac{1}{8760533070643200\,000},\text{ }b(5,5)=\tfrac{227}{%
12854317559387145633792000000}, \\
b(1,6) &=&\tfrac{1}{127288\,050\,516627\,176\,816640\,000\,000}, \\
b(6,6) &=&\tfrac{133933}{25516459094444104187401241999966208000000000}, \\
b(1,7) &=&\tfrac{1}{998707926079695101611943783301120000000000}, \\
b(7,7) &=&\tfrac{2006509}{895370835179\,281010%
\,419215815294340559070476369920000000000000}.
\end{eqnarray*}

Having an explicit formula for the $b(h,k)$ and $b(k,k)$ would via (\ref{kli}%
) help to decide whether the conjecture $\lambda =\infty $ is true subject
to the Riemann hypothesis. Using (\ref{kli}) and the values of $b(1,k)$ and $%
b(k/k)$, we have the new lower values $\Lambda (k)$ for $k=2,...,7.$

\begin{equation*}
\begin{tabular}[t]{|l|l|l|l|l|l|}
\hline
$\Lambda (2)$ & $\Lambda (3)$ & $\Lambda (4)$ & $\Lambda (5)$ & $\Lambda (6)$
& $\Lambda (7)$ \\ \hline
$1.375\,3$ & $1.8858$ & $2.343\,9$ & $2.764\,0$ & $3.149\,1$ & $3.5004$ \\ 
\hline
\end{tabular}%
\end{equation*}

Next in the following, we will apply the Wirtinger inequality due to Agarwal
and Pang \cite{AP} to establish a new explicit formula for the lower bounds
of $\Lambda .$ This inequality is presented in the following theorem.

\bigskip

\textbf{Theorem C}. \textit{Assume that }$x(t)\in C^{1}[0,$\textit{\ }$\pi ]$%
\textit{\ and }$x(0)=x(\pi )=0$\textit{, then }%
\begin{equation}
\int_{0}^{\pi }(x^{^{\prime }}(t))^{2k}dt\geq \frac{2\Gamma \left(
2k+1\right) }{\pi ^{2k}\Gamma ^{2}\left( \left( 2k+1\right) /2\right) }%
\int_{0}^{\pi }x^{2k}(t)dt\text{, \ for }k\geq 1,  \label{AP}
\end{equation}

\textbf{Theorem 2.3.} \textit{Assuming the Riemann hypothesis and the moment
(\ref{Hu}) is correctly predicted we have }%
\begin{equation}
\Lambda (k)\geq \frac{1}{2\pi }\left( \frac{b(0,k)}{b(k,k)}\frac{2\Gamma
\left( 2k+1\right) }{\Gamma ^{2}\left( \left( 2k+1\right) /2\right) }\right)
^{\frac{1}{2k}},\text{ \ for }k=3,4,...\text{ }.  \label{APo}
\end{equation}%
\textbf{Proof. }To prove this theorem we will employ\textbf{\ }the
inequality (\ref{AP}). By a suitable linear transformation, we can deduce
from (\ref{AP}) that: if $x(t)\in C^{1}[a,b]$ and $x(a)=x(b)=0,$ then 
\begin{equation}
\int_{a}^{b}(\frac{b-a}{\pi })^{2k}(x^{^{\prime }}(t))^{2k}dt\geq \frac{%
2\Gamma \left( 2k+1\right) }{\pi ^{2k}\Gamma ^{2}\left( \left( 2k+1\right)
/2\right) }\int_{a}^{b}x^{2k}(t)dt\text{, \ for }k\geq 1.  \label{AP0}
\end{equation}%
Since by our assumption (\ref{Hu}) is correctly predicted by RMT, we have
for $k=h,$ the moments of the derivative of $Z(t)$%
\begin{equation}
\int_{0}^{T}(Z^{^{\prime }}(t))^{2k}dt\sim a(k)b(k,k)T\left( \log T\right)
^{k^{2}+2k},  \label{Z1}
\end{equation}%
and for $h=0,$ we have the moments of $Z(t)$ 
\begin{equation}
\int_{0}^{T}Z^{2k}(t)dt\sim a(k)b(0,k)T\left( \log T\right) ^{k^{2}}.
\label{Z2}
\end{equation}%
Now, we follow the proof of \cite{Hall3} by supposing that $t_{l}$ is the
first zero of $Z(t)$ not less than $T$ and $t_{m}$ the last zero not greater
than $2T.$ Suppose further that for $l\leq n<m$, we have 
\begin{equation}
L_{n}=t_{n+1}-t_{n}\leq \frac{2\pi \kappa }{\log T},  \label{R131}
\end{equation}%
and applying the inequality (\ref{AP0}), to obtain 
\begin{equation*}
\int_{t_{n}}^{t_{n+1}}\left[ \left( \frac{L_{n}}{\pi }\right) ^{2k}\left(
Z^{^{\prime }}(t)\right) ^{2k}-\frac{2\Gamma \left( 2k+1\right) }{\pi
^{2k}\Gamma ^{2}\left( \left( 2k+1\right) /2\right) }Z^{2k}(t)\right] dt\geq
0.
\end{equation*}%
Since the inequality remains true if we replace $L_{n}/\pi $ by $2\kappa
/\log T$, we have 
\begin{equation}
\int_{t_{n}}^{t_{n+1}}\left[ \left( \frac{2\kappa }{\log T}\right)
^{2k}\left( Z^{^{\prime }}(t)\right) ^{2k}-\frac{2\Gamma \left( 2k+1\right) 
}{\pi ^{2k}\Gamma ^{2}\left( \left( 2k+1\right) /2\right) }Z^{2k}(t)\right]
dt\geq 0.  \label{RR1}
\end{equation}%
Summing (\ref{RR1}) over $n,$ applying (\ref{Z1}) and (\ref{Z2}), we obtain 
\begin{eqnarray*}
&&a(k)b(k,k)\left( \frac{2\kappa }{\log T}\right) ^{2k}T\left( \log T\right)
^{k^{2}+2k}-\frac{2a(k)b(0,k)\Gamma \left( 2k+1\right) }{\pi ^{2k}\Gamma
^{2}\left( \frac{2k+1}{2}\right) }T\left( \log T\right) ^{k^{2}} \\
&=&\left( a(k)b(k,k)\kappa ^{2k}(2^{2k})-\frac{2a(k)b(0,k)\Gamma \left(
2k+1\right) }{\pi ^{2k}\Gamma ^{2}\left( \left( 2k+1\right) /2\right) }%
\right) T\left( \log T\right) ^{k^{2}} \\
&\geq &O(T\log ^{k^{2}}T),
\end{eqnarray*}%
whence, as $T\rightarrow \infty $, we obtain 
\begin{equation*}
\kappa ^{2k}\geq \frac{a(k)b(0,k)}{2^{2k}a(k)b(k,k)}\frac{2\Gamma \left(
2k+1\right) }{\pi ^{2k}\Gamma ^{2}\left( \frac{2k+1}{2}\right) }=\frac{b(0,k)%
}{2^{2k}b(k,k)}\frac{2\Gamma \left( 2k+1\right) }{\pi ^{2k}\Gamma ^{2}\left(
\left( 2k+1\right) /2\right) }.
\end{equation*}%
This implies that 
\begin{equation*}
\Lambda ^{2k}(k)\geq \frac{b(0,k)}{2^{2k}b(k,k)}\frac{2\Gamma \left(
2k+1\right) }{\pi ^{2k}\Gamma ^{2}\left( \left( 2k+1\right) /2\right) },
\end{equation*}%
which is the desired inequality. The proof is complete.

\bigskip

Having an explicit formula for the $b(k,k)$ would via (\ref{APo}) help to
decide whether the conjecture $\lambda =\infty $ is true subject to the
Riemann hypothesis. Using (\ref{APo}) and the values of $b(0,k)/b(k/k)$ (see
Table 2), we have the new lower values of $\Lambda $%
\begin{equation}
\begin{tabular}[t]{|l|l|l|l|l|}
\hline
$\Lambda (3)$ & $\Lambda (4)$ & $\Lambda (5)$ & $\Lambda (6)$ & $\Lambda (7)$
\\ \hline
$2.2265$ & $2.6544$ & $3.0545$ & $3.4259$ & $3.7676$ \\ \hline
\end{tabular}%
.  \label{SA1}
\end{equation}

In the following, we will apply the Wirtinger inequality due Brneti\'{c} and
Pe\v{c}ari\'{c} \cite{BP} to establish a new explicit formula for the lower
bounds of $\Lambda .$

\bigskip

\textbf{Theorem 2.4.} \textit{Assuming the Riemann hypothesis and the moment
(\ref{Hu}) is correctly predicted, we have}%
\begin{equation}
\Lambda (k)\geq \frac{1}{2\pi }\left( \frac{b(0,k)}{b(k,k)}\frac{1}{I(k)}%
\right) ^{\frac{1}{2k}}\text{, \ \ for }k=3,4,...\text{ .}  \label{SA2}
\end{equation}

\textbf{Proof.} To prove this theorem, we will employ the inequality (\ref%
{AP2}). Proceeding as in Theorem 2.2, we may have%
\begin{equation*}
\kappa ^{2k}\geq \frac{a(k)b(0,k)}{2^{2k}a(k)b(k,k)}\frac{1}{\pi ^{2k}I(k)}=%
\frac{b(0,k)}{2^{2k}b(k,k)}\frac{1}{\pi ^{2k}I(k)},\text{ (as }T\rightarrow
\infty ).
\end{equation*}%
This implies that 
\begin{equation*}
\Lambda ^{2k}(k)\geq \frac{b(0,k)}{2^{2k}b(k,k)}\frac{1}{\pi ^{2k}I(k)},
\end{equation*}%
which is the desired inequality. The proof is complete.

\bigskip

Again having an explicit formula for the $b(k,k)$ would via (\ref{SA2}) help
to decide whether the conjecture $\lambda =\infty $ is true subject to the
Riemann hypothesis. To find the new estimation of $\Lambda (k)$ we need the
values of $I(k)$ for $k=3,...,7$, which are calculated numerically in the
following table:%
\begin{equation*}
\begin{tabular}[t]{|l|l|l|l|l|}
\hline
$I(3)$ & $I(4)$ & $I(5)$ & $I(6)$ & $I(7)$ \\ \hline
$\frac{19\,581}{5000\,000}$ & $\frac{743}{1000\,000}$ & $\frac{14\,961}{%
100\,000\,000}$ & $\frac{15\,653}{500\,000\,000}$ & $\frac{16\,823}{%
2500\,000\,000}$ \\ \hline
\end{tabular}%
\end{equation*}%
Using these values and the values in Table 2, we have by using the explicit
formula (\ref{SA2}) the new estimation of $\Lambda (k)$ in the following
table:%
\begin{equation*}
\begin{tabular}[t]{|l|l|l|l|l|}
\hline
$\Lambda (3)$ & $\Lambda (4)$ & $\Lambda (5)$ & $\Lambda (6)$ & $\Lambda (7)$
\\ \hline
$2.4905$ & $2.9389$ & $3.350\,8$ & $3.728\,7$ & $4.0736.$ \\ \hline
\end{tabular}%
\end{equation*}

\end{document}